\documentclass[11pt]{article}
\usepackage{amsfonts}
\usepackage{amsmath}
\usepackage{amssymb}
\begin{document}

\baselineskip 18pt
\title{\Large\bf On power residues modulo a prime}
\author{Ke Gong and Chaohua Jia}
\date{}
\maketitle {\small\noindent {\bf Abstract.} Let $p$ be a
sufficiently large prime number, $n$ be a positive odd integer with
$n|\,p-1$ and $n>p^\varepsilon $, where $\varepsilon$ is a
sufficiently small constant. Let $k(p,\,n)$ denote the least
positive integer $k$ such that for
$x=-k,\,\dots,\,-1,\,1,\,2,\,\dots,\,k$, the numbers $x^n\pmod p$
yield all the non-zero $n$-th power residues modulo $p$. In this
paper, we shall prove
$$
k(p,\,n)=O(p^{1-\delta}),
$$
which improves a result of S. Chowla and H. London in the case of
large $n$.}

\vskip.3in
\noindent{\bf 1. Introduction}

Let $p>3$ be a prime number, and $n$ be a positive odd integer with
$n|\, p-1$. Let $k(p,\,n)$ denote the least positive integer $k$
such that for $x=-k,\,\dots,\,-1,\,1,\,2,\,\dots,\,k$, the numbers
$x^n\pmod p$ yield all the non-zero $n$-th power residues modulo
$p$.

Obviously, one has
$$
k(p,\,n)<\frac12 p.
$$
S. Chowla and H. London [2], using elementary arguments, proved that
$$
{p-1\over 2n}\leq k(p,\,n)<({1\over 2}-{1\over 2n})p.
$$

In this paper, we improve the result of S. Chowla and H. London in
the case of large $n$.

\medskip\noindent
{\bf Theorem}. \emph{Let $p$ be a sufficiently large prime number,
$n$ be a positive odd integer with $n|\,p-1$ and $n>p^\varepsilon $,
where $\varepsilon $ is a sufficiently small constant. The number
$k(p,\,n)$ is defined as above. Then we have}
$$
k(p,\,n)=O(p^{1-\delta}),
$$
\emph{where $\delta=\delta(\varepsilon)$ is a sufficiently small
constant which depends only on $\varepsilon$}.

\vskip.3in
\noindent{\bf 2. The proof of Theorem}

\noindent {\bf Lemma 1}. \emph{Assume that $(m,\,p)=1$, $x_0$ is a
solution of the congruence equation}
$$
x^n\equiv m\pmod p. \eqno (1)
$$
\emph{Then the equation (1) has $n$ different solutions modulo $p$
in the following form:}
$$
x_0 g^{{p-1\over n}j},\qquad\quad j=0,\,1,\,\dots,\,n-1,
$$
\emph{where $g$ is a primitive root modulo $p$.}

For a reference, see [2].

Let $\mathbb{F}_p$ be the finite field with $p$ elements,
$\mathbb{F}_p^\ast=\mathbb{F}_p\setminus\{0\}$,
$$
H=\{h:\ h\equiv g^{{p-1\over n}j}\pmod p,\ j=0,\,1,\,\dots,\,n-1\}.
$$
Obviously $H$ is a multiplicative subgroup of $\mathbb{F}_p^\ast$.
The order of $H$ is $n$. The numbers $x_0h$ $(h\in H)$ yield all the
solutions of the equation (1), where $x_0$ is one of the solutions
of the equation (1).

\medskip\noindent
{\bf Lemma 2}~(Bourgain--Glibichuk--Konyagin [1]). \emph{For any
$\varepsilon>0$ and for any multiplicative subgroup $H$ of
$\mathbb{F}_p^\ast$ with $|H|>p^\varepsilon$, there exists a small
positive constant $\delta=\delta(\varepsilon)$ such that}
$$
\max_{a\in\mathbb{F}_p^\ast}\,\Bigl|\sum_{h\in H}
e\Bigl({ah\over p}\Bigr)\Bigr| \ll{|H|\over p^{3\delta}},
$$
where $e(x)=e^{2\pi ix}.$

\medskip\noindent
\textit{The proof of Theorem.} We have
\begin{align*}
&\ \,\sum_{h\in H}\sum_{\substack{|x|\leq p^{1-\delta}\\ x\ne 0\\
x_0h\equiv x\,({\rm mod}\,p)}} 1 \\
& = \sum_{h\in H}\sum_{\substack{|x|\leq p^{1-\delta}\\ x\ne
0}}{1\over p}\sum_{r=1}^p e\Bigl({r(x_0h-x)\over p}\Bigr) \\
& = {1\over p}\sum_{r=1}^p \Bigl(\sum_{h\in H}e\Bigl({rx_0 h\over
p}\Bigr)\Bigr)\Bigl(\sum_{\substack{|x|\leq p^{1-\delta}\\ x\ne
0}}e\Bigl(-{rx\over p}\Bigr)\Bigr)\\
& = {|H|\over p}\cdot 2[p^{1-\delta}]+{1\over p}\sum_{r=1}^{p-1}
\Bigl(\sum_{h\in H}e\Bigl({rx_0 h\over p}\Bigr)\Bigr)\Bigl(
\sum_{\substack{|x|\leq p^{1-\delta}\\ x\ne 0}}e\Bigl(-{rx\over
p}\Bigr)\Bigr).
\end{align*}
The application of Lemma 2 yields
\begin{align*}
&\ \,{1\over p}\sum_{r=1}^{p-1}
\Bigl(\sum_{h\in H} e\Bigl({rx_0 h\over p}\Bigr)\Bigr)
\Bigl(\sum_{\substack{|x|\leq p^{1-\delta}\\ x\ne 0}}
e\Bigl(-{rx\over p}\Bigr)\Bigr) \\
& \ll {1\over p}\sum_{r=1}^{p-1}\Bigl|\sum_{h\in H}e\Bigl({rx_0
h\over p}\Bigr)\Bigr|\Bigl|\sum_{\substack{|x|\leq p^{1-\delta}\\
x\ne 0}}e\Bigl(-{rx\over p}\Bigr)\Bigr| \\
& \ll {|H|\over p^{1+3\delta}}\sum_{r=1}^{p-1}\Bigl|\sum_{\substack{
|x|\leq p^{1-\delta}\\ x\ne 0}}e\Bigl(-{rx\over p}\Bigr)\Bigr| \\
& \ll {n\over p^{1+3\delta}}\sum_{r=1}^{p-1}{1\over \|{r\over p}\|} \\
& \ll {n\over p^{1+3\delta}}\sum_{1\leq r\leq{{p-1}\over 2}}{p\over r} \\
& \ll {n\over p^{2\delta}},
\end{align*}
where $\|\alpha\|$ denotes the distance from $\alpha$ to the nearest integer.

It follows that
$$
\sum_{h\in H}\sum_{\substack{|x|\leq p^{1-\delta}\\ x\ne 0\\
x_0h\equiv x\,({\rm mod}\,p)}} 1
= {2n\over p^\delta}+O\Bigl({n\over p^{2\delta}}\Bigr)
\geq {n\over p^\delta}
\geq 1.
$$
Hence there is a solution $x_0h$ of the equation (1) such that
$$
x_0h\equiv x\,({\rm mod}\,p),\qquad\quad |x|\leq p^{1-\delta},
$$
which implies
$$
k(p,\,n)\leq p^{1-\delta}.
$$

So far the proof of Theorem is complete.

\vskip.3in
\noindent{\bf Acknowledgements}

Ke Gong is supported by the National Natural Science Foundation of
China (Grant No. 11671119). Chaohua Jia is supported by the National
Natural Science Foundation of China (Grant No. 11771424 and Grant
No. 11321101).

\vskip.3in
\noindent{\bf Additional remark}

After the first version of this paper was published in arXiv, we
found that the result in Theorem of this paper had been included in
Theorem 8 of the paper of M. Ram Murty (\emph{Small solutions of
polynomial congruences}, Indian J. Pure Appl. Math., \textbf{41}
(2010), no.1, 15-23).

\vskip.6in

\noindent
Ke Gong \\
Department of Mathematics, Henan University, Kaifeng, Henan 475004,
P. R. China \\
E-mail: {\tt kg@henu.edu.cn}

\

\noindent
Chaohua Jia \\
Institute of Mathematics, Academia Sinica, Beijing 100190, P. R.
China \\
Hua Loo-Keng Key Laboratory of Mathematics, Chinese Academy of
Sciences, Beijing 100190, P. R. China \\
School of Mathematical Sciences, University of Chinese Academy of
Sciences, Beijng 100049, P. R. China \\
E-mail: {\tt jiach@math.ac.cn}


\vskip.6in

\begin{thebibliography}{9}

\bibitem{1} J. Bourgain, A. A. Glibichuk and
S. V. Konyagin, \emph{Estimates for the number of sums and products
and for exponential sums in fields of prime order}, J. London Math.
Soc. \textbf{73} (2006), 380--398.

\bibitem{2} S. Chowla and H. London, \emph{Bounds on
the $n$-th power residues $\pmod p$}, Canad. Math. Bull. \textbf{12}
(1969), 679--680.


\end{thebibliography}
\end{document}